On Harrod Neutral Technical Change.
Uzawa Revisited.

Thomas Russell
Santa Clara University
Dept. of Economics

Dec 2004

There has recently been a renewed interest in the theorem of Uzawa (1961) that in the classic Solow/Swan growth model, only Harrod neutral (labor augmenting) technical progress is consistent with an equilibrium steady state. Jones and Scrimgeour (2004) for example, show that a constant rate of labor augmenting technical change is a consequence of the economy being in dynamic equilibrium, there being no need to assume (as in Barro and Sala I Martin (1995)) that the underlying rate of technical change is already constant.

What is still perhaps not completely clear is just why technical change is so constrained. After all, to a mathematician, the fact that output $Y$ can be written as
$$Y = f(K, L, t)$$
(where, as usual, $K$ is capital, $L$ is labor, and $t$ captures the effect of technical progress), seems to permit a high degree of interchangeability between the factors of production, particularly when we assume, as we do, that the production function is homogeneous of degree one in capital and labor. It may seem a little odd that equilibrium growth forces this function to have the form $Y = f(K, L^{e^{\lambda t}})$, the Harrod form.

In this note we give a two line proof of this proposition which is designed to separate the role played by a priori assumptions concerning the shape of the function f, assumptions which hold both in and out of the steady state, from properties which hold because the economy is indeed in long run dynamic steady state.

Assumptions: We will make all the standard neoclassical assumptions set out in Jones and Scrimgeour (op cit.). In particular labor $L$ will be assumed to grow at the constant rate n. The proof then follows from the development of the Solow Swan model in Sato(1964). That is to say, given the production relationship
$$Y(t) = f(K(t), L(t), t)$$
the dynamics must satisfy
$$\frac{\dot{Y}}{Y} = \frac{f_K K}{Y} \frac{\dot{K}}{K} + \frac{f_L L}{Y} \frac{\dot{L}}{L} + \frac{f_t}{Y} \qquad (1)$$
this by total differentiation of Y, see Sato (op.cit.) p 381.

We now make two assumptions, one with respect to the mathematical form of the function $f$, the other with respect to the long run economic behavior of the system.

**Assumption 1**: The function $f$ is homogeneous of degree 1 so that by Euler's[1] theorem $\frac{f_K K}{Y} = 1 - \frac{f_L L}{Y}$.

**Assumption 2:** The long run behavior of output is characterized by a steady state in which output and capital both grow at the same constant rate g.

---

[1] Credit for noticing the importance of Euler's theorem for the neoclassical theory of the distribution of income belongs to Flux (1894). This long forgotten paper can still be read with much pleasure for its late Victorian style as much as its extremely important content.

Then we have

**Theorem:** In the steady state output must satisfy the equation $f_t + nf_L(t)L(t) = g$, where $g$ is a constant. This is a standard partial differential equation of mathematical physics (the inhomogeneous variable coefficient advective equation) whose solution is easily found by the method of characteristics. Thus so we have

**Corollary** (Uzawa) Output in steady state growth is given by $Y=f(K(t), L(t) \exp^{(g-n)t})$, i.e. output grows at the constant rate $g = n + \rho$ where $\rho$ is the constant rate of labor augmenting (Harrod neutral) technical progress.

Proof.: The proof is immediate. By homogeneity $\frac{f_K K}{Y} = 1 - \frac{f_L L}{Y}$. Then the fact that we are in the steady state means that,

$$\frac{\dot{Y}}{Y} = \frac{\dot{K}}{K} = n + \frac{f_t}{f_L L} = g, \text{ some constant.} \quad (2)$$

Mathematicians recognize the last two terms in (2) as the inhomogeneous form of the variable coefficient advective equation $f_t = \alpha f_L(t)L(t) = 0$ noted above. The solution is known[2] to be $Y=f(K(t), L(t)\exp^{(g-n)t})$

**Conclusion:** What is the economic import of this theorem? Presumably any interest in properties of the steady state must be associated with the belief that an economy with an arbitrary initial condition finds itself close to the steady state in some reasonable amount of time. Unfortunately, as we know from the work of Sato(1964), Atkinson (1969), and others, for reasonable parameter values, the time to reach the steady state can exceed the lifetime of all of the all those alive at any starting time. In these circumstances the value of the Uzawa theorem is not what it has to say about the behavior of any real economy, but rather what is implies about the shape of the transition paths along which we actually live. We have examined techniques for studying these transition paths elsewhere, Russell and Zecevic (1998).

Also it should be noted that technical change in the Solow/Swan model is exogenous. If we make technical change endogenous, then it is possible to have biased technical change with constant factor shares even in the steady state, see Sato and Ramachandran (1999) and Sato, Ramachandran, and Lian (2000)

---

[2] Sarra (2003)